\documentclass{amsart}
\usepackage[latin1]{inputenc}
\usepackage{latexsym}
\usepackage{amsfonts}
\usepackage{amssymb}
\usepackage{color}
\usepackage{dsfont}
\usepackage{varioref}
\usepackage{epsfig}
\usepackage{xy}
\usepackage[active]{srcltx}
\usepackage{array}
\RequirePackage{graphics,epsfig,psfrag}
\usepackage{color}
\usepackage{colordvi,fancyhdr}
\usepackage{pdfpages}
\definecolor{Gris}{cmyk}{0.1,0.1,0.1,.75}
\usepackage[colorlinks=true,citecolor=Gris,linkcolor=Gris,filecolor=Gris,urlcolor=Gris]{hyperref}
\def\imagetop#1{\vtop{\null\hbox{#1}}}
\def\pn{\medskip\par\noindent}

\def\frac#1#2{{\textstyle{{#1} \overwithdelims.. {#2}}}}
\def\Frac#1#2{{\displaystyle{{#1} \overwithdelims.. {#2}}}}
\def\abs#1{\left \vert #1 \right \vert}
\def\eps{\varepsilon}
\def\Sum{\mathop{\sum}\limits}
\newcommand{\Pf}{\noindent{\em Proof}. }

\newcommand{\EPf}{\hbox{}\hfill$\Box$\vspace{.5cm}}

\newtheorem{thm}{Theorem}[section]
\newtheorem{defi}[thm]{Definition}
\newtheorem{prop}[thm]{Proposition}
\newtheorem{lemma}[thm]{Lemma}
\newtheorem{cor}[thm]{Corollary}
{\theoremstyle{definition}
\newtheorem{rem}[thm]{Remark}}{\theoremstyle{definition}
\newtheorem{exa}[thm]{Example}}

\newcommand{\ZZ}{{\mathbb Z}}

\def\Mod#1{\,(\hbox{\rm mod}\,#1)}
\def\[#1\]{\begin{equation} #1 \end{equation}}

\begin{document}
\title{Untangling trigonal diagrams}
\date{\today}
\author{Erwan Brugall\'e}
\author{Pierre-Vincent Koseleff}
\author{Daniel Pecker}

\begin{abstract}
Let
$K$
be a link
of Conway's normal form $C(m)$, $m \geq 0$, or $C(m,n)$ with $mn>0$,
and let
$D$ be a trigonal diagram of $K.$
We show that
it is possible to transform
$D$ into an alternating trigonal diagram, so that
all intermediate
diagrams remain trigonal, and the number of crossings never increases.
\end{abstract}
\maketitle

\section{Introduction}
\begin{figure}[!ht]
\begin{center}
\begin{tabular}{ccc}
{\scalebox{.8}{\includegraphics{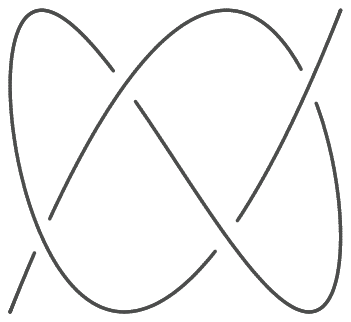}}}&\quad&
{\scalebox{.8}{\includegraphics{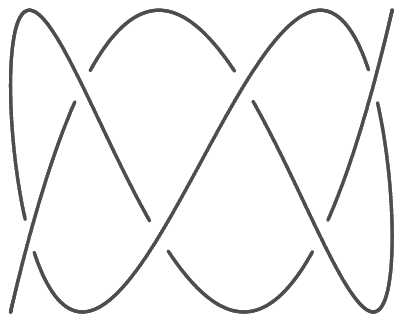}}}\\
$4_1=C(2,2)$&&$5_1=C(5)$
\end{tabular}
\end{center}
\caption{Polynomial representations  of the knots $4_1$ and $5_1$}
\label{ctrex}
\end{figure}

If we try to simplify a knot
or link diagram,
then
the number of crossings may have to be increased in some intermediate diagrams,
see \cite{G,KL,A,Cr}.
In this paper, we shall see that this
strange
phenomenon
cannot occur
for trigonal diagrams of two-bridge torus links
and for generalized twist links.
The next theorem is the main result of this paper.
\begin{thm} \label{thm:isotopy}
Let $K$ be a link of Conway's normal form $C(m)$, $m \geq 0$, or $C(m,n)$ with $mn>0$,
and let $D$ be a trigonal diagram of $K$.
Then, it is possible to transform $D$ into an alternating trigonal diagram,
so that all intermediate diagrams remain trigonal, and the number of crossings never increases.
\end{thm}
We also prove that if $K$ is a two-bridge link which is not of these two types,
then $K$ admits diagrams that cannot be simplified without increasing the number of crossings.
Our original motivation to tackle this problem is the study of polynomial knots,
 their polynomial isotopies and their degrees,
see examples in Figure \ref{ctrex}, \cite{BKP2,KP1,KP3,RS,V}.
In
particular this is why we
prefer to consider our knots as long knots. As an application
 of Theorem \ref{thm:isotopy}, we determine in \cite{BKP2} the
lexicographic degree of  two-bridge knots of Conway's normal form $C(m)$
with $m$ odd, or $C(m,n)$ with $mn$ positive and even.

\medskip
The paper is organized as follows.
In Section \ref{sec:2-bridge}, we recall Conway's notation for trigonal diagrams
of two-bridge links, and their classification by their Schubert fractions.
In Section \ref{sec:simple} we  define {\em slide isotopies} as
trigonal isotopies such that the number of crossings never increases.
We find necessary conditions for a  two-bridge link diagram to be  {\em simple},
that is to say it
cannot
be transformed into a simpler diagram by any slide isotopy.
We use continued fraction properties to prove
Theorem \ref{thm:isotopy} in Section \ref{sec:proof}.
In Section \ref{sec:awkward} we
show that if a two-bridge link is neither a torus link  nor a twist link,
then it possesses awkward trigonal diagrams.

\section{Trigonal diagrams  of two-bridge knots}\label{sec:2-bridge}
\begin{figure}[!ht]
 \centerline{
 \psfrag{a}{$m-1$}\psfrag{b}{$1-n$}
 {\scalebox{.70}{\includegraphics{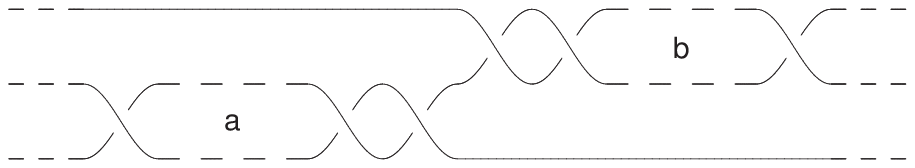}}} \quad
 \psfrag{a}{$m-1$}\psfrag{b}{$n-1$}
 {\scalebox{.70}{\includegraphics{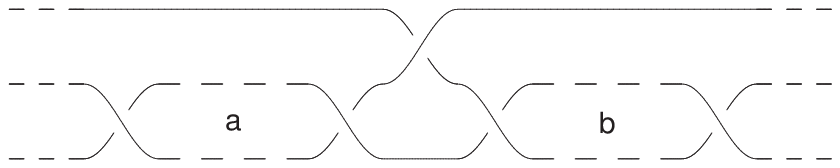}}}
}
 \caption{Lagrange isotopies}
 \label{lagrange}
 \end{figure}

A two-bridge link  admits a diagram in \emph{Conway's  open form} (or trigonal form).
This diagram, denoted by
$D(m_1, m_2, \ldots, m_k)$  where $m_i  \in \ZZ^*$ are integers,
is explained by the following picture (see \cite{Co}, \cite[p.~187]{Mu}).
\psfrag{a}{\small $m_1$}\psfrag{b}{\small $m_2$}
\psfrag{c}{\small $m_{k-1}$}\psfrag{d}{\small $m_{k}$}
\begin{figure}[!ht]
\begin{center}
{\scalebox{.75}{\includegraphics{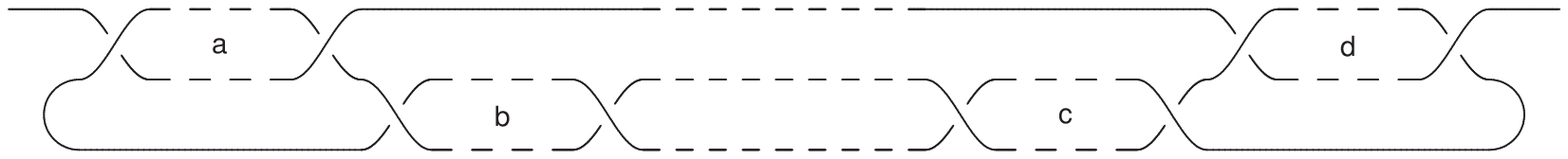}}}\\[30pt]
{\scalebox{.75}{\includegraphics{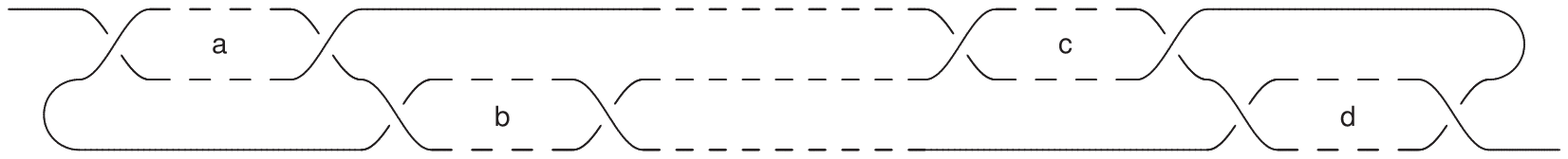}}}
\end{center}
\caption{Conway's  form for long
links}
\label{conways3}
\end{figure}
The number of twists is denoted by the integer $\abs{m_i}$, and the sign of $m_i$ is defined
as follows: if $i$ is odd, then the right twist is positive,
if $i$ is even, then the right twist is negative.
In  Figure \ref{conways3} the $m_i$ are positive (the $m_1$ first twists are right twists).
\smallskip
These diagrams are also called $3$-strand-braid representations, see \cite {KL1,KL}.
\pn
The two-bridge links are classified by their Schubert fractions
$$
\Frac {\alpha}{\beta} =
m_1 + \Frac{1} {m_2 + \Frac{1} {\cdots +\Frac 1{m_k}}}=
[ m_1, \ldots, m_k], \quad \alpha >0, \ (\alpha,\beta)=1.
$$
Given $[ m_1, \ldots, m_k]=\Frac {\alpha}{\beta} $ and
$[ m'_1, \ldots, m'_l]=\Frac {\alpha'}{\beta'}$, the diagrams
 $D(m_1, m_2, \ldots, m_k)$ and $D(m'_1, m'_2, \ldots, m'_l)$
correspond to isotopic links if and only if
$\alpha = \alpha' $ and $ \beta' \equiv \beta ^{\pm 1} \Mod{\alpha}$,
see \cite[Theorem 9.3.3]{Mu}.
The integer $ \alpha$ is odd for a knot, and even for a two-component link.
\pn
Any fraction admits a continued fraction expansion
$
\Frac {\alpha}{\beta} = [m_1, \ldots, m_k]
$
where all the  $m_i's $ have the same sign.
Therefore every two-bridge link $K$
admits a diagram in \emph{Conway's normal form}, that is an alternating
diagram of the form
$ D( m_1, m_2, \ldots m_k) $ where all the  $m_i's $ have the same
sign. In this case we will write $L= C(m_1, \ldots , m_k).$
\pn
It is classical that one can transform any trigonal diagram of a two-bridge link into its
Conway's normal form using the  Lagrange isotopies,  see \cite{KL} or
\cite[p.~ 204]{Cr}:
\[
D( x,m,-n,-y) \rightarrow  D(x,m-\eps,\eps,n-\eps,y), \ \eps=\pm 1,
\]
where $m,n$ are integers, and $x,y$  are  sequences of
integers  (possibly empty), see Figure \ref{lagrange}.
 These isotopies twist  a part of the diagram, and  the number of crossings may increase in  intermediate
diagrams.
Since we want to simplify links without  increasing their complexity, we introduce
different isotopies in the following section.


\section{Slide isotopies, simple and awkward  diagrams}\label{sec:simple}
\begin{figure}[!ht]
\begin{center}
{\scalebox{.6}{\includegraphics{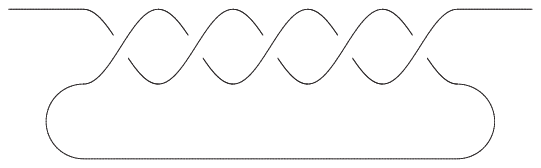}}} \quad
{\scalebox{.6}{\includegraphics{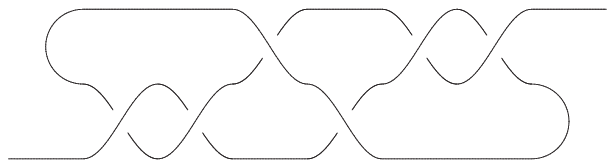}}} \quad
{\scalebox{.6}{\includegraphics{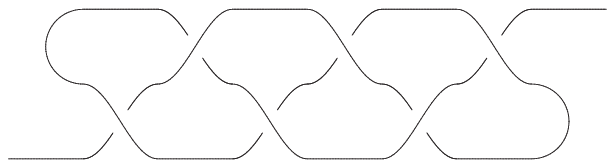}}} \\[20pt]
{\scalebox{.6}{\includegraphics{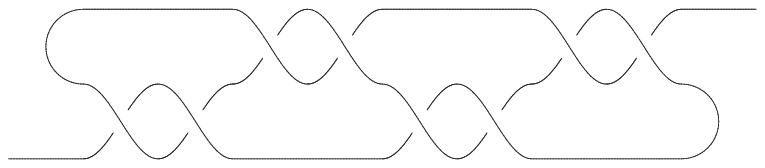}}} \quad
{\scalebox{.6}{\includegraphics{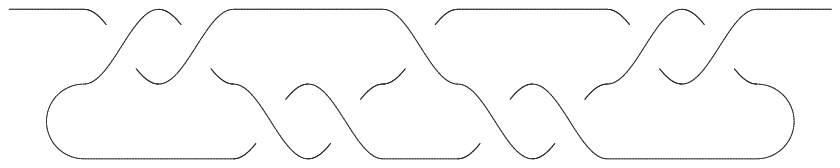}}}
\end{center}
\caption{Some trigonal diagrams of the torus knot $5_1$}
\label{fig:51}
\end{figure}

\begin{defi}
We shall say that an isotopy of trigonal diagrams is a {\em slide} isotopy
if the number of crossings never increases and if
all the intermediate diagrams remain trigonal.
\end{defi}
\begin{exa}Some diagrams of the torus knot $5_1 = C(5)$:
$D(5)$, $D(2,1,-1,-2)$, $D(-1,-1,-1,1,1,1)$,
$D(-2,2,-2,2)$,  and $ D(2, 2, -1, 2,2)$ are depicted in Figure \ref{fig:51}.
By  Theorem \ref{thm:isotopy}, it is possible to simplify  these diagrams
into the alternating diagram $D(5)$  by slide isotopies only.
\end{exa}
\begin{rem}(G\"{o}ritz, 1934)
Let $D$ be a trigonal diagram of a link $K$.
Let $\widetilde{K}$ be the image of $K$
by a half-turn around the $x$-axis,
and let $\widetilde{D}$ be the $xy$-projection of $\widetilde{K}$, see Figure \ref{half-turn}.
\begin{figure}[!th]
\begin{center}
{\scalebox{.75}{\includegraphics{conway2.eps}}}\\[30pt]
{\scalebox{.75}{\includegraphics{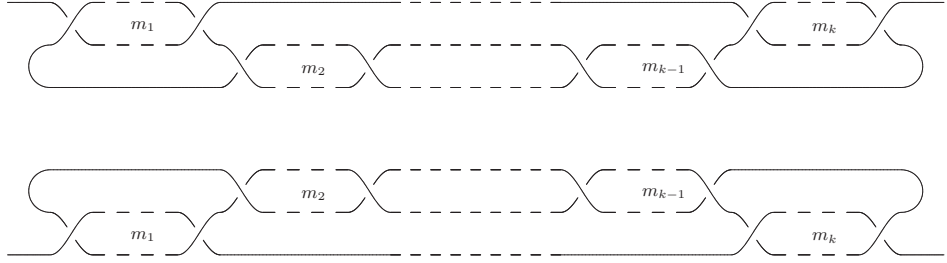}}}
\end{center}
\caption{The two
diagrams $D$ and $\widetilde{D}$  with the same Conway notation.}
\label{half-turn}
\end{figure}
The diagrams $D$ and $\widetilde{D}$ are  diagrams of the same link $K,$
nevertheless they are generally not isotopic by a slide isotopy (\cite{G}).
It is often convenient to identify the diagrams $D$ and $\widetilde{D}.$
\end{rem}
\begin{defi}
We define the complexity of a trigonal diagram $D(m_1, \ldots, m_k)$ as
$c(D)=k+ \Sum \abs{m_i}.$
\end{defi}
\begin{defi} A trigonal diagram is called a {\em simple} diagram if it cannot
be simplified into a diagram
of lower complexity by using slide isotopies only. A non-alternating simple diagram
is called an {\em awkward} diagram.
\end{defi}
The next example shows the existence of awkward diagrams.
\begin{exa}
Let us consider the diagram $D= D(4,-3).$ It is an awkward  diagram of the
knot $6_2=C(3,1,2)$: the only possible Reidemeister moves increase the
number of crossings.
Of course, we can transform  this diagram into an alternating one
using  Lagrange isotopies, but in this process some intermediate diagrams
will have more crossings than $D$.
\end{exa}
\begin{figure}[!ht]
\begin{center}
\begin{tabular}{ccc}
{\scalebox{.8}{\includegraphics{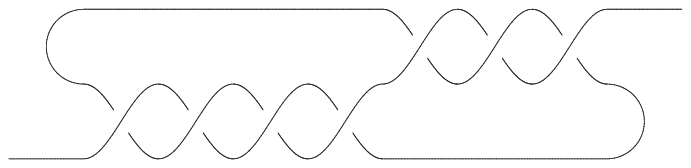}}} &\quad&
{\scalebox{.8}{\includegraphics{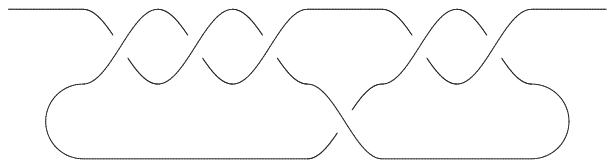}}}\\
$D(4,-3)$ & & $6_2 = C(3,1,2)$
\end{tabular}
\end{center}
\caption{$D(4,-3)$ is an awkward  diagram of the knot $6_2 = C(3,1,2)$ }
\end{figure}
\pn
\begin{rem}\label{rem:awkward}
More generally, let $m_1, \ldots m_k$ be   integers
that are neither all positive nor all negative, and such that $ |m_i| \ge 2$,
$ \left ( |m_1| \ge 3 \right .$ or $ \left . m_1m_2>0 \right ) $ and
$ \left ( |m_k| \ge 3 \right .$ or $\left .  m_{k-1}m_k >0 \right ).$
Then the diagram
$D(m_1, \ldots , m_k) $ is awkward. In fact the only Reidemeister moves that can be applied increase
the number of crossings.
Kauffman and Lambropoulou call such diagrams {\em hard  diagrams}, see \cite{KL}.
\end{rem}
\begin{prop}\label{prop:simple}
Let $ D=D(m_1, \ldots , m_k) , \ k >1$ be a simple trigonal  diagram.
Then we have
\begin{enumerate}
\item[$(i)$]
$\abs{m_1} \ge 2$;  $  \abs{m_k} \ge 2$;
${m_i} \ne 0,  \  i=2,\ldots, k-1$;
\item[$(ii)$] $m_1 m_2 >0$ or $\abs{m_1}\ge 3$; $ \  m_{k-1}m_k>0$ or $\abs{m_k} \ge 3$;
\item[$(iii)$] for $i = 2,\ldots,k$,
$m_{i-1} m_i \ne -1$;
\item [$ (iv)$]
suppose that $\abs {m_i}= 1$, then
\begin{enumerate}
\item if $m_{i-1} m_i <0$, then
$i\le k-2$,  $m_i m_{i+1} > 0,$ and $m_i m_{i+2} > 0$;
\item if $m_{i} m_{i+1} <0$, then
$i\ge 3$,   $m_{i-2} m_{i} > 0$, and $m_{i-1} m_{i} > 0$.
\end{enumerate}
\end{enumerate}
\end{prop}
\Pf
\par
\begin{figure}[!ht]
\begin{center}
\psfrag{m}{$m$}\psfrag{m+1}{$m+1$}
{\scalebox{.8}{\includegraphics{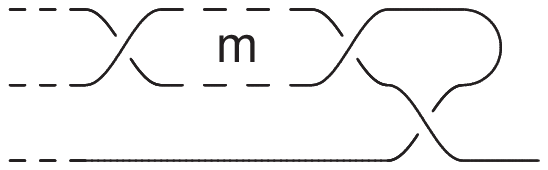}}} \quad
{\scalebox{.8}{\includegraphics{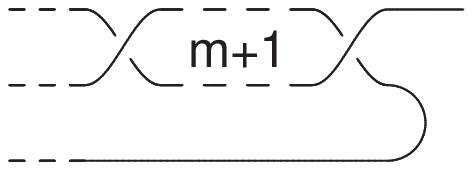}}}
\end{center}
\caption{The slide isotopy $D(x, m,1) \rightarrow D( x, m+1)$}
\label{fig:(i)}
\end{figure}

\begin{figure}[!ht]
\begin{center}
\psfrag{a}{$m-1$}
\centerline{
{\scalebox{.8}{\includegraphics{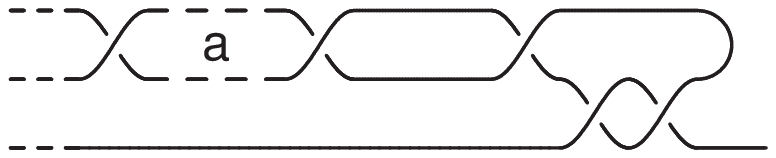}}}\quad
{\scalebox{.8}{\includegraphics{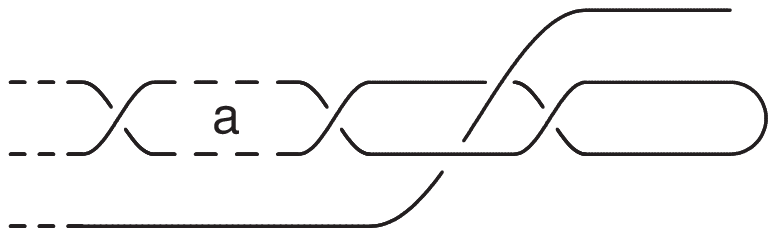}}}}
\hbox{}\hfill
{\scalebox{.8}{\includegraphics{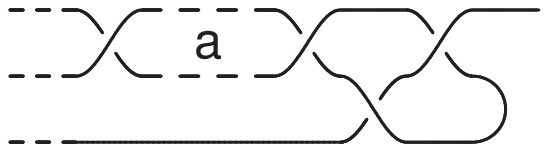}}}
\end{center}
\caption{The slide isotopy $ D (x, m, -2) \rightarrow D(x, m-1,2),  \ m>0$}
\label{fig:iso2}
\end{figure}

\begin{figure}[!ht]
\begin{center}
\psfrag{a}{$m$}\psfrag{b}{$n$}\psfrag{c}{$p$}
{\scalebox{.8}{\includegraphics{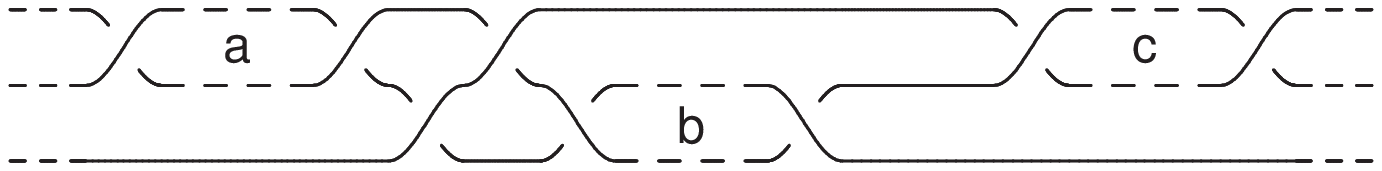}}} \\[10pt]
\psfrag{a}{$m$}\psfrag{b}{$-n$}\psfrag{c}{$p$}
{\scalebox{.8}{\includegraphics{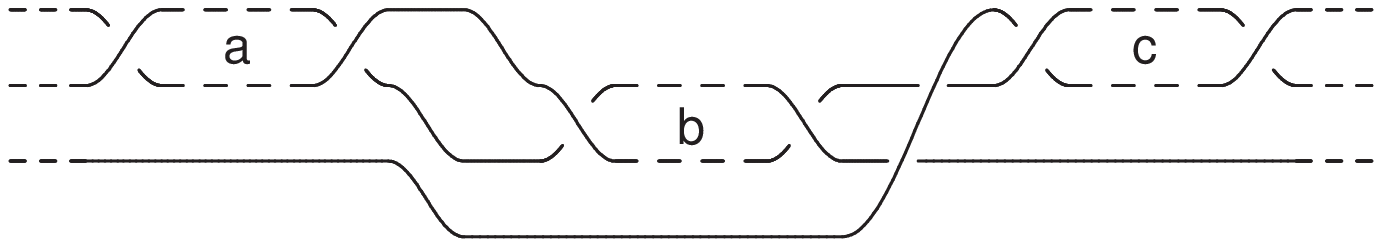}}} \\[10pt]
\psfrag{a}{$m-n$}\psfrag{b}{$p+1$}
{\scalebox{.8}{\includegraphics{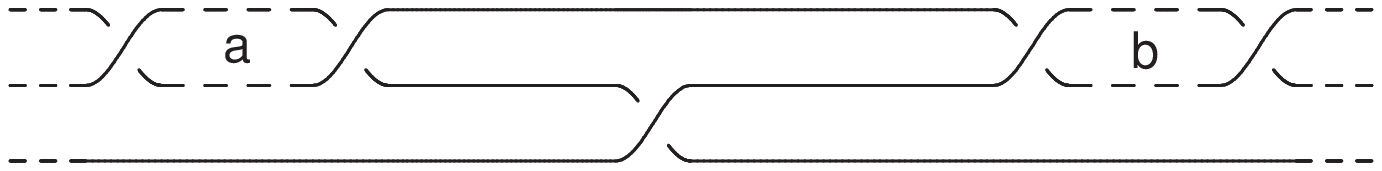}}}
\end{center}
\caption{The slide isotopy  $ D( x,m,-1,1,n,p,y) \rightarrow D ( x,m-n,-1, 1+p, y)$}
\label{fig:iso3}
\end{figure}

\begin{figure}[th]
\begin{center}
\begin{tabular}{cc}
\psfrag{a}{$m-1$}\psfrag{b}{$n$}
{\scalebox{.6}{\includegraphics{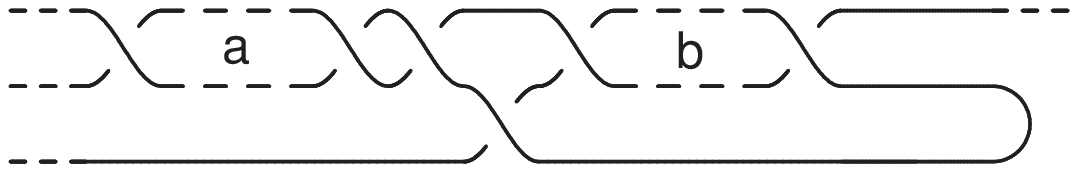}}} &
\psfrag{a}{$m-1$}\psfrag{b}{$-n$}
{\scalebox{.6}{\includegraphics{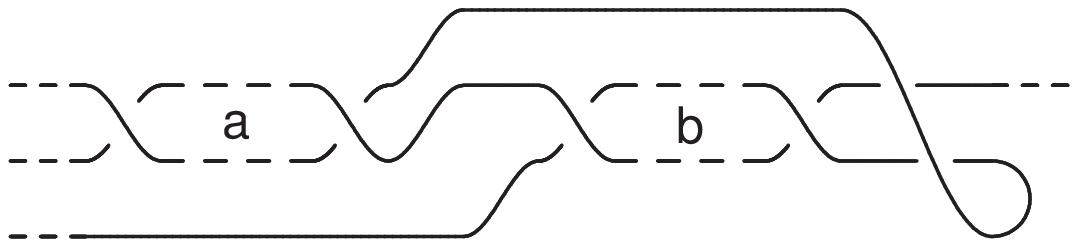}}} \\[20pt]
\psfrag{a}{$m-1$}\psfrag{b}{$-n+1$}
{\scalebox{.6}{\includegraphics{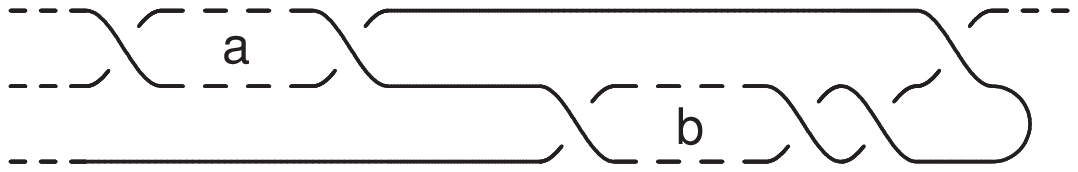}}} &
\psfrag{a}{$m-1$}\psfrag{b}{$-n+1$}
{\scalebox{.6}{\includegraphics{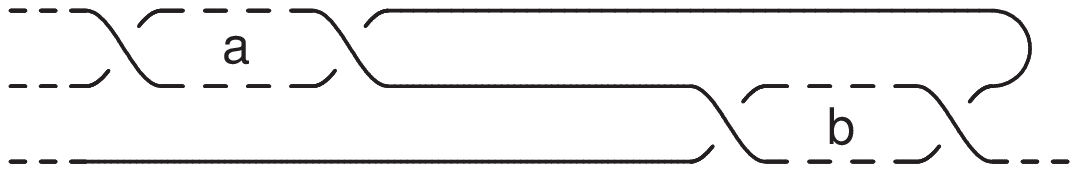}}}
\end{tabular}
\end{center}
\caption{The slide isotopy $D(x,m,-1,n) \mapsto D(x,m-1,-n+1)$, $m>0$}
\label{fig:iso4}
\end{figure}

\begin{figure}[!ht]
\psfrag{a}{$m-1$}\psfrag{b}{$n$}\psfrag{c}{$p$}
{\scalebox{.8}{\includegraphics{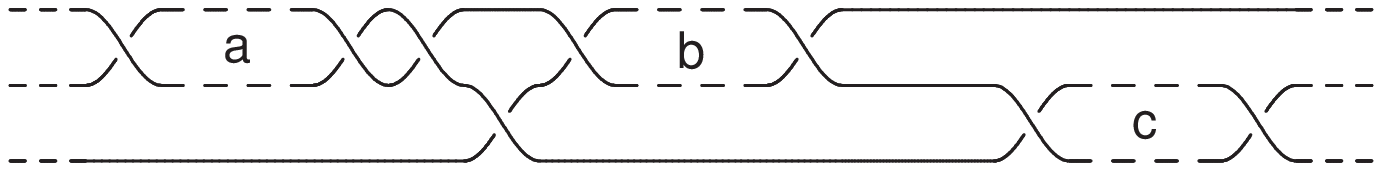}}} \\[20pt]
\psfrag{a}{$m-1$}\psfrag{b}{$-n$}\psfrag{c}{$p$}
{\scalebox{.8}{\includegraphics{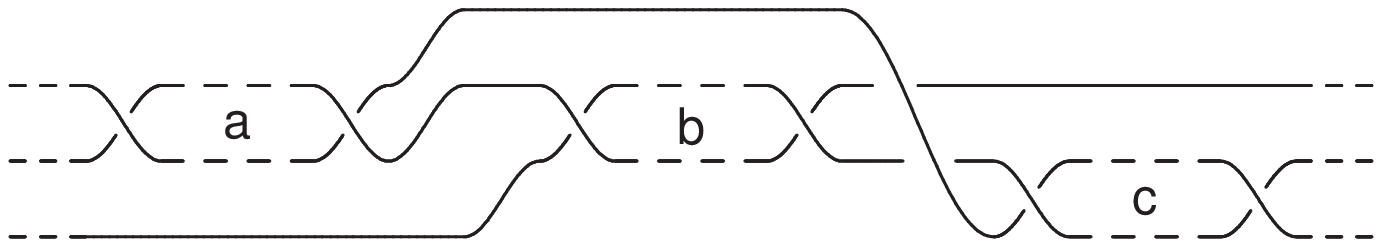}}} \\[20pt]
\psfrag{a}{$m-1$}\psfrag{b}{$-n$}\psfrag{c}{$p$}
{\scalebox{.8}{\includegraphics{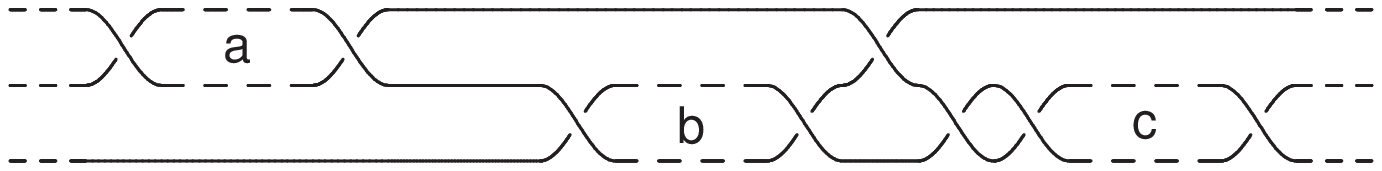}}}
\caption{The slide isotopy  $ D( x,m,-1,n,p,y) \rightarrow D ( x,m-1,-n,1,p-1, y)$, $m>0$, $p<0$}
\label{fig:iso5}
\end{figure}

\noindent
$(i)$
the slide isotopy $ D(x, m, 0) \rightarrow D(x)$ diminishes the complexity,
consequently $m_k \ne 0,$ and similarly $m_1 \ne 0.$
The slide isotopy $D(x, m,1) \rightarrow D( x, m+1)$ diminishes the
complexity, see Figure \ref{fig:(i)}.
Consequently
$m_k \ne 1$ since $D$ is simple.
We also have $m_k \ne -1, $ and similarly $m_1 \ne \pm1.$
The slide isotopy
$D( x, m, 0, n, y) \rightarrow D (x,m+n, y) $ shows that $m_i \ne 0.$
\pn
$(ii)$
Let us  show that if $|m_k|=2$ then $m_k m_{k-1} >0$.
Suppose on the contrary that, for example, $m_k=-2$ and $m_{k-1}>0$. Then the
slide isotopy $ D (x, m, -2) \rightarrow D(x, m-1,2)$
 decreases the complexity of $D$, see Figure \ref{fig:iso2}.
This contradicts the simplicity of $D$.
Similarly, we see that if $|m_1|=2$ then $m_1 m_2 >0.$
\pn
$(iii)$
Consider the slide isotopy
$ D( x,m,-1,1,n,p,y) \rightarrow D ( x,m-n,-1, 1+p, y),$ see Figure \ref{fig:iso3}.
When $(p,y)=\emptyset$, this isotopy becomes $D (x,m,-1,1,n) \rightarrow D(x, m-n-1).$
\def\imagetop#1{\vtop{\null\hbox{#1}}}
\def\imagebot#1{\vbot{\null\hbox{#1}}}
It lowers the complexity, and consequently a simple diagram
cannot be of the form $ D(u,-1,1,v).$
\pn
$(iv)$
By $(i)$ we have  $i\ne 1$ and $i\ne k$.
Let us assume that $m_{i-1}m_i<0$,
the proof in the  case $m_im_{i+1}<0$ being entirely
similar. The  slide isotopy
$ D (x, m, -1, n ) \rightarrow  D( x, m-1, -n+1) $ depicted in  Figure \ref{fig:iso4}
shows that $i\le k-2$.
Since $D$ is minimal,
the slide isotopy
$ D(x, m , -1, n, p, y) \rightarrow  D ( x, m-1, -n, 1, p-1, y ) $
depicted in Figure \ref{fig:iso5} implies that $p<0,$ that is
$ m_im_{i+2}>0$.

Now, assume  that $m_im_{i+1}<0$.
Let $i$ be  the maximal integer  such that there exists a simple diagram $ D(m_1, \ldots m_k)$ of length $k$
such that $ \abs{m_i} =1,  \  m_{i-1} m_i <0, $ and $ m_im_{i+1} <0.$
Once more we use the slide isotopy
$ D(x, m , -1, n, p, y) \rightarrow  D ( x, m-1, -n, 1, p-1, y ) $
depicted in Figure \ref{fig:iso5}.
Since $p= m_{i+2} >0$ and  $n= m_{i+1} >0$ by assumption, the new diagram contradicts the maximality
of  $i.$
\EPf
\pn
The condition $(iv)$ of Proposition \ref{prop:simple}
asserts that a simple diagram
cannot contain any subsequence
$ \pm ( m, -1, n)$, $m,n >0$.
Nevertheless, it can contain subsequences of the form $\pm (m,-1)$, $m>1$.
However, the next corollary shows that this phenomenon can be avoided.
\begin{cor} \label{cor:isotopy}
Let $D$ be a trigonal Conway diagram of a
two-bridge link.
Then it is possible to transform $D$ by slide isotopies into a simple
diagram $ D( m_1, \ldots , m_k) $ such that
for $i=2, \ldots ,k$, either  $\abs{ m_i} \ne 1, \ $  or  $ m_{i-1} \, m_i  >0$.
\end{cor}
\Pf
Let us suppose that each simple diagram deduced from $D$
by slide isotopies contains subsequences of the form
$\pm (m,-1), \  m \ge 2,$
and let $\mu \ne 0$  be the minimum number of such subsequences.

Among
all simple diagrams $D(m_1, \ldots , m_k)$ deduced from $D$ by slide isotopies and possessing
 $\mu$ such subsequences
there is a minimal integer $r$ such that
 $ (m_{r-1}, m_r) = \pm ( m, -1), \ m \ge 2$.
By Proposition \ref{prop:simple},
$(iv)$ we have $r\leq  k-2$ and our diagram is of the form
$\Delta=D( x, m, -1, n, p, y)$, where $n,p<0$.
The use of the slide isotopy
depicted in Figure \ref{fig:iso5}:
$\Delta \rightarrow \Delta'= D( x, m-1, -n, 1, p-1, y )$,
contradicts the minimality of either $\mu$ or $r$.
\EPf
\section{Proof of
Theorem \ref{thm:isotopy}}\label{sec:proof}
The proof is based on some arithmetical properties of the continued fractions
$[m_1, \ldots , m_k]$ related to
simple diagrams $D(m_1, \ldots, m_k)$, that are
consequences of Proposition \ref{prop:simple}
and
Corollary
\ref{cor:isotopy}.
\begin{lemma}\label{lem:schubert}
Let $x$ be a rational number
defined by its continued fraction expansion
$x= \Frac {\alpha}{\beta} = [ m_1, m_2, \ldots , m_k]$,
$m_1>0$, $|m_k| \ge 2$, $m_i \in {\bf Z}^{*}$, $\alpha \ge 0$.
Suppose that for $i = 2 , \ldots, k$,
we have either $| m_i| \ne 1 $  or $m_{i-1} m_i > 0.$
Then the following hold:
\begin{itemize}
\item[$(a)$] $x> m_1-1,$ and consequently $\alpha >0$  and $ \beta >0$;
\item[$(b)$] if $k \ge 2$ and
  ($m_{k-1} m_k > 0$ or $\abs{m_k} \ge 3$), then
  $ \alpha \ge 2$ and $\beta \ge 2$;
\item[$(c)$]
if
in addition we
have $m_1 \ge 2$ and
$( m_1 m_2 > 0$ or $m_1 \ge 3 )$ then $ x>2$;
\item[$ (d)$]
if in addition we
have $k \ge 3$ and $( \abs{m_2} \ne 1 $ or $m_2m_3 > 0), $  then
$ \alpha \not \equiv 1  \Mod{\beta}.$
\end{itemize}
\end{lemma}
\Pf
\pn
$(a)$ We use an induction on $k$.
If $k=1,$ then we have $x= m_1 > m_1-1,$ and the result is true.
Let us suppose that $k \ge 2.$
If $m_2 >0,$ then
 we have $ y= [m_2, \ldots , m_k] >0,$ and then
$ x = m_1 +1/ y > m_1 >m_1-1.$
If $m_2<0,$ then by assumption $-m_2 \ge 2$. By induction we have
$ y= [ -m_2, \ldots , -m_k] >1,$ and then
$x=m_1-1/ y > m_1-1.$
\pn
$(b)$
We use again an induction on $k$.
If $k=2$, then there are two cases to consider.

If $m_2 >0,$  then we have $ \beta= m_2 \ge 2,$ and $ \alpha= m_1m_2+1 \ge 2.$

If $m_2<0, $ then we have $ \abs{m_2} \ge 3.$
Then $ \beta= |m_2| \ge 3$ and $ \alpha= m_1 |m_2|-1 \ge 3 m_1-1 \ge 2.$
\pn
Let us suppose now that $k \ge 3.$
Let $x=m_1 + q/p,$ where $p/q= [m_2, \ldots , m_k], \ p>0.$

If $m_2 >0$, then $\beta=p \ge 2$ and $q\ge 2$ by induction,
and then  $ \alpha = m_1p + q
\ge 2.$

If $m_2<0,$ then $-m_2 \ge 2.$
By induction, $ [ -m_2-1, -m_3, \ldots , -m_k ] = \Frac {p+q}{-q} $
is such that $ p+q \ge 2$ and $-q \ge 2 .$
Consequently,  we obtain
$ \beta= p
\ge 2-q
\ge 4 \ge 2, $ and $ \alpha= m_1 p + q  = m_1 (p+q) + (m_1-1) (-q) \ge 2m_1 \ge 2.$
\pn
$(c)$
If $m_2>0,$ then $ x > m_1 \ge 2.$
If $m_2<0,$ then $m_1 \ge 3,$ and then $ x> m_1-1 \ge 2$ by $(a).$
\pn
$ (d)$
If $m_2 <0,$ then $-m_2 \ge 2$ and
as in the proof of $(b)$ we get  $p+q \ge 2 $ and $-q\ge 2$, where
 $ \Frac p{-q} = [-m_2, \ldots , -m_k] >1.$
Consequently we obtain $ p \ge 2-q>1-q>0$    and then $ q \not \equiv 1  \  \Mod{p} .$
As $ \alpha = m_1p+q$ and $\beta= p,$ the result is proved in this case.

Suppose now that $ m_2>0$, and consider $ \Frac u{v}= [m_3, \ldots ,
  m_k]$ with $u>0$. In particular, we have $ \beta = m_2 u +v$, and  $\alpha= m_1
\beta + u \equiv u \   \Mod{\beta}$.

If $m_3>0$, then $v>0$ and so  $ \beta >u$.
If $k=3,$  then $u= m_3 \ge 2.$ If $k>3, $ then we  have $u \ge 2 $ by $(b).$
Hence
$ \beta >u \ge 2 ,$
and so
$\alpha \equiv u \not \equiv 1   \   \Mod{\beta} .$

If
$m_3 <0,$ then $m_2 \ge 2$, $ -m_3 \ge 2 $, and $v<0$.
We have
$\Frac u {-v}>1$
by $(a),$ i.e. $u>-v$.
Hence we have $ \beta \ge 2u +v > u \ge 2,$
and we obtain $\alpha \equiv u \not \equiv 1  \    \Mod{\beta}.$
\EPf
\pn
{\em Proof of Theorem \ref{thm:isotopy}}.
Let $D(m_1, \ldots, m_k) $ be a simple diagram deduced from
$D$
by
admissible isotopies, and
satisfying the condition of
Corollary \ref{cor:isotopy}. Recall that it also satisfies all conclusions of
Proposition \ref{prop:simple}.
We write $\Frac{\alpha}{\beta} = [m_1,m_2, \ldots , m_k]$ with
$\alpha>0$ and $(\alpha,\beta) = 1$.
Without loss of generality, we may assume that both $m$ and $n$ are
nonnegative.

\pn
Let us first consider the case when  $K$ is the torus link
$C(m)$. By the classification of torus links by their Schubert fractions, we
have $\alpha=m$ and $\beta =1\Mod m$. If $k\ge 2$, then
Lemma \ref{lem:schubert} $(b)$ implies that $|\beta|\ge 2$ so $\abs \beta \geq m-1$.
But then, by Lemma \ref{lem:schubert} $(c)$, we would have $\alpha = m > 2 |\beta|\geq 2m-2$  so $m\leq 1$.
Hence we obtain $k=1$, which proves the result.
\pn
Now, let us consider the case when  $K$ is the twist link
$C(m,n)$ with $m\ge2$ and $n\ge 2$.
Then we have $\alpha=mn+1$, and either
$\beta \equiv n\Mod\alpha$ or $\beta \equiv -m\Mod\alpha$.
Lemma \ref{lem:schubert} $(c)$ implies that $\alpha > 2|\beta|$, from
which we deduce that $\beta = n$ or $\beta = -m$.
Consequently we have $\alpha \equiv 1 \Mod{\beta}$, which implies by
Lemma \ref{lem:schubert} $(d)$ that $k=2$.

By Proposition \ref{prop:simple} $(i)$, we have $|m_1|\ge 2$ and  $|m_2|\ge 2$.
If  $m_1m_2 < 0,$ then we would have
$ \Frac{\alpha}{\beta}= \Frac{ |m_1m_2| -1 } { \pm m_2}$. In
particular we would have
$\alpha \equiv -1 \Mod{\beta},$ which is impossible since $\alpha
\equiv 1 \Mod{\beta}$
 and
$ |\beta | = |m_2| \ge 3$ by Proposition \ref{prop:simple} $(ii)$.
Consequently we have $m_1m_2 >0$ and the
diagram is alternating.
\EPf
\section{Some awkward trigonal diagrams}\label{sec:awkward}
\begin{figure}[!ht]
\psfrag{a}{$m$}\psfrag{b}{$1-a$}\psfrag{c}{$-n$}
{\scalebox{.6}{\includegraphics{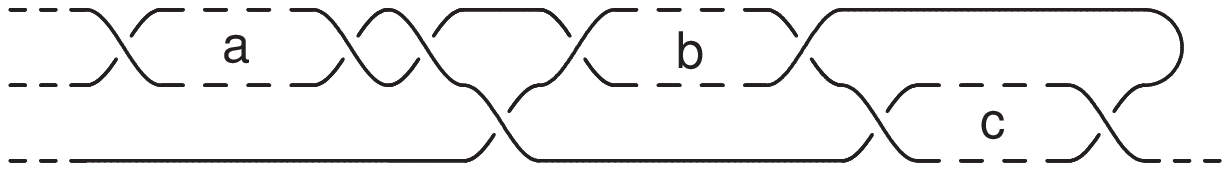}}} \\[20pt]
\psfrag{a}{$m$}\psfrag{b}{$1-a$}\psfrag{c}{$-n$}
{\scalebox{.6}{\includegraphics{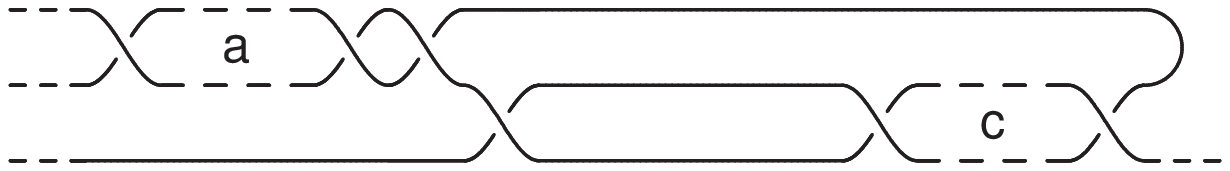}}} \\[20pt]
\caption{The awkward diagrams $D(x, m+1,-1,1-a,-n)$ and
$D(x, m+1,-n-1)$.}
\label{fig:awkward}
\end{figure}

The following  result  shows  that if a two-bridge
link
is not of the Conway normal form  $C(m)$, or
$
C(m,n)$ with $mn>0$, then it possesses an awkward trigonal diagram.
\begin{prop}
Let $k \ge 3$ and let $K$ be a two-bridge
link of Conway normal form $C(m_1, \ldots , m_k)$, $m_i >0$, $m_1 \ge 2, $ $m_k \ge 2.$
Then $K$ possesses an awkward trigonal diagram.
\end{prop}
\Pf
Let $[m_{k-2},m_{k-1},m_k]=[m,a,n]$.
Using the Lagrange identity we have
$[m,a,n]=[m+1,-1,1-a,-n]$. If $a=1$, then this last continued fraction
is $[m+1,-n-1]$.
Therefore $K$ admits the trigonal diagram
\begin{itemize}
\item[] $D(m_1, \ldots , m_{k-3}, m+1,-1,1-a,-n)$, if $a>1$;
\item[or ] $D(m_1, \ldots , m_{k-3}, m+1,-n-1)$ if $a=1$.
\end{itemize}
These two diagrams are awkward (see Figure \ref{fig:awkward}).
\EPf
\pn
Furthermore, we have a stronger result.
\begin{thm}\label{thm:awkward}
Let $k \ge 3$ and let $K$ be a two-bridge
link of Conway normal form $C(m_1, \ldots , m_k)$, $m_i >0$, $m_1 \ge 2, $ $m_k \ge 2.$
Then $K$ possesses a hard trigonal diagram.
\end{thm}
\Pf
Using the identity
$\Frac{az+1}{(a-1)z+1} =
[\underbrace{2,-2,2, \ldots, (-1)^{a-2}2}_{a-1 \hbox{ \scriptsize terms}},(-1)^{a-1}(z+1)]$, we obtain
$$
[m,a,n] = [m+1,-\Frac{an +1}{(a-1)n+1}] =
[m+1,-2,2,-2, \ldots, (-1)^{a-1}2,(-1)^a(n+1)],
$$
and we deduce that
$$[m_1, \ldots, m_k] = [x,m,a,n] =
[x,m+1, -2, 2, \ldots , (-1)^{a-1} 2, (-1)^{a} (n+1)].$$
Therefore, the diagram $D = D(x, m+1, -2, 2, \ldots , (-1)^{a-1} 2, (-1)^{a} (n+1))$
is a diagram of $K$, and it is a hard diagram  by  Remark \ref{rem:awkward}.
\EPf
\begin{rem}
Theorem \ref{thm:isotopy}
shows that the trivial knot has no trigonal awkward diagram.
On the other hand, G\"{o}ritz  found an awkward diagram of the trivial
knot in \cite {G}, and Kauffman and Lambropoulou constructed many such examples
(see  \cite {KL, A, Cr}).
\end{rem}

\goodbreak
\vfill
\pn
\hrule width 7cm height 1pt 
\pn
{\small
Erwan {\sc Brugallé}\\
École Polytechnique,
Centre Mathématiques Laurent Schwartz, 91 128 Palaiseau Cedex, France\\
e-mail: {\tt erwan.brugalle@math.cnrs.fr}
\pn
Pierre-Vincent {\sc Koseleff}\\
Universit\'e Pierre et Marie Curie (UPMC Sorbonne Universit\'es),\\
Institut de Math\'ematiques de Jussieu (IMJ-PRG)  \& Inria-Rocquencourt\\
e-mail: {\tt koseleff@math.jussieu.fr}
\pn
Daniel {\sc Pecker}\\
Universit\'e Pierre et Marie Curie (UPMC Sorbonne Universit\'es),\\
Institut de Math\'ematiques de Jussieu (IMJ-PRG),\\
e-mail: {\tt pecker@math.jussieu.fr}
}

\end{document}